%&amstex          
\input amstex\documentstyle{amsppt}  
\pagewidth{12.5cm}\pageheight{19cm}\magnification\magstep1
\topmatter
\title Parabolic character sheaves and Hecke algebras\endtitle
\author G. Lusztig\endauthor
\address{Department of Mathematics, M.I.T., Cambridge, MA 02139}\endaddress
\thanks{Supported by NSF grant DMS-2153741}\endthanks
\endtopmatter   
\document

\define\pos{\text{\rm pos}}
\define\Irr{\text{\rm Irr}}

\define\frl{\forall}
\define\pe{\perp}
\define\si{\sim}

\define\sqc{\sqcup}

\define\bA{\bar A}

\define\bZ{\bar Z}
\define\lb{\linebreak}

\define\op{\oplus}
   
\define\part{\partial}
\define\em{\emptyset}

\define\m{\mapsto}
\define\do{\dots}

\define\lra{\leftrightarrow}

\define\sub{\subset}    
\define\bxt{\boxtimes}
\define\T{\times}
\define\ti{\tilde}
\define\nl{\newline}
\redefine\i{^{-1}}
\define\fra{\frac}
\define\un{\underline}

\define\ot{\otimes}
\define\bbq{\bar{\QQ}_l}

\define\End{\text{\rm End}}

\define\tr{\text{\rm tr}}

\redefine\spa{\spadesuit}

\define\a{\alpha}

\redefine\c{\chi}
\define\g{\gamma}
\redefine\d{\delta}
\define\e{\epsilon}

\define\p{\pi}
\define\ph{\phi}

\define\s{\sigma}
\redefine\t{\tau}

\redefine\l{\lambda}

\define\x{\xi}

\define\vt{\vartheta}

\define\kk{\bold k}

\define\NN{\bold N}

\define\PP{\bold P}
\define\QQ{\bold Q}

\define\VV{\bold V}

\define\ZZ{\bold Z}

\define\YY{\bold Y}

\define\ca{\Cal A}
\define\cb{\Cal B}

\define\cd{\Cal D}
\define\ce{\Cal E}

\define\ch{\Cal H}

\define\ck{\Cal K}
\define\cl{\Cal L}
\define\cm{\Cal M}

\define\cp{\Cal P}

\define\ct{\Cal T}

\define\cv{\Cal V}
\define\cw{\Cal W}

\define\fC{\frak C}

\define\tc{\ti c}

\define\tE{\ti E}

\define\tG{\ti G}

\define\tV{\ti V}

\define\tZ{\ti Z}

\define\sha{\sharp}

\define\bg{\bar g}

\head Introduction\endhead
\subhead 0.1\endsubhead
Let $\kk$ be an algebraic closure of a finite field $F_q$ with $q$
elements.
Let $G$ be a reductive algebraic group over $\kk$ with a fixed
$F_q$-structure. In \cite{L76,L77} it has been shown
that the endomorphism
algebra of a representation of $G(F_q)$ induced from a unipotent
cuspidal representation of the reductive quotient of a parabolic
subgroup defined over $F_q$ is an Iwahori-Hecke algebra $\ch_q$ with
possibly unequal parameters; moreover the parameters were explicitly
described. In \cite{L03,\S27} a new realization of $\ch_q$ was given
in terms involving functions on a variety which was later used to
define parabolic character sheaves \cite{L04}. This new realization
was refined in \cite{L06, 36.16} in terms of actual parabolic
character sheaves at least in the case where $G$ is a symplectic
group, but the statements there were given without proof. In
this paper a proof of those statements (in a slightly modified form)
is given. We will show elsewhere that a similar approach applies to
any $G$ (assumed to be split over $F_q$).

\subhead 0.2\endsubhead
The affine analogues of the algebras $\ch_q$ were introduced in
\cite{L91,\S17}. The methods of this paper give a realization of
these algebras in terms of the affine parabolic sheaves studied in
\cite{L10}.

\subhead 0.3. Notation\endsubhead
For any connected affine algebraic group $H$ we denote by $U_H$ the
unipotent radical of $H$ and by $R_H$ the solvable radical of $H$.

We set $\ca=\ZZ[v,v\i]$ where $v$ is an indeterminate. For any
$\ca$-module $\cm$ we write $\un\cm=\QQ(v)\ot_\ca\cm$.

If $j:X@>>>Y$ is a smooth morphism of algebraic varieties with
connected fibres of dimension $d$ we set $j^\star(A)=j^*(A)[d]$ for any
complex $A$ of $\bbq$-shaves on $Y$.

If $F:X@>>>X$ is a map of sets we write $X^F=\{x\in X;F(x)=x\}$.

For $n\in\NN$ we set $[0,n]=\{z\in\NN;z\le n\}$.

\head 1. Recollection on parabolic character sheaves\endhead
\subhead 1.1\endsubhead
Let $\tG$ be a (possibly disconnected) reductive group over $\kk$ with
identity component $G$ and with a fixed connected component $D$.
Let $W$ be the Weyl group of $G$ and let $I$ be the set of simple
reflections of $W$. Let $l:W@>>>\NN$ be the standard length function.
Let $\t:W@>>>W$ be the automorphism induced by the automorphism of $G$
given by conjugation by an element of $D$. We have $\t(I)=I$.

For any $J\sub I$ let $W_J$ be the subgroup of $W$ generated by $J$.
For $J\sub I,K\sub I$ let

${}^JW=\{w\in W;l(yw)=l(y)+l(w){} \frl y\in W_J\}$,

$W^K=\{w\in W;l(wy)=l(y)+l(w){} \frl y\in W_K\}$,

${}^JW^K={}^JW\cap W^K$.

Let $\cp$ be the set of parabolic subgroups of $G$.
Let $\cb\sub\cp$ be the set of Borel subgroups of $G$.
Now $W$ can be regarded as the set of orbits of $G$ acting by
simultaneous conjugation on $\cb\T\cb$. For $(B,B')\in\cb\T\cb$ let
$\pos(B,B')$ be the Weyl group element corresponding to the orbit of
$(B,B')$. To $P\in\cp$ we associate the subset $J$ of $I$ consisting
of all  $i\in I$ such that for some $(B,B')\in\cb\T\cb$ with
$B\sub P$, $B'\sub P$ we have $\pos(B,B')=i$. For $J\sub I$ let
$\cp_J$ be the set of all $P\in\cp$ with associate subset $J$. We have
$\cp=\sqc_{J\sub I}\cp_J$, $\cb=\cp_\em$.
For $J\sub I,K\sub I$ and $(P,Q)\in\cp_J\T\cp_K$ there is a unique
element $w=\pos(P,Q)\in{}^JW^K$ such that $\pos(B,B')\in W_JwW_K$
for any $(B,B')\in\cb\T\cb$ with $B\sub P,B'\sub Q$.

Let $J\sub I$. For $w\in W$ let $w_J$ be the unique element of
minimal length in $W_{\t(J)}wW_J$. We have a partition
${}^{\t(J)}W={}^{\t(J)}_*W\sqc{}^{\t(J)}_\spa W$ where
$${}^{\t(J)}_*W=\{w\in{}^{\t(J)}W;w_JJw_J\i=\t(J)\},$$
$${}^{\t(J)}_\spa W=\{w\in{}^{\t(J)}W;w_JJw_J\i\ne\t(J)\}.$$
Using results of B\'edard in \cite{B85},\cite{L04,\S2}, we see that

(a) if $w\in{}^{\t(J)}_*W$, then $w=w_J$ hence $wJw\i=\t(J)$.

\subhead 1.2\endsubhead
For $P\in\cp,Q\in\cp$ we set $P^Q=(P\cap Q)U_P$; we have $P^Q\in\cp$,
$P^Q\sub P$ hence $U_P\sub U_{P^Q}$.

\subhead 1.3\endsubhead
Following \cite{L04}, for $J\sub I$ we set
$$Z_J=\{(P,P',gU_P);(P,P')\in\cp_J\T\cp_{\t(J)},gU_P\in D/U_P,
gPg\i=P'\}.$$
For $w\in W$ let
$$\tZ_\em^w=\{(B,B',g)\in\cb\T\cb\T D;gBg\i=B',\pos(B',B)=w\}.$$
Define $\p^w:\tZ_\em^w@>>>Z_J$ by $(B,B',g)\m(P,P',gU_P)$ where
$P\in\cp_J,P'\in\cp_{\t(J)}$, $B\sub P,B'\sub P'$. A simple perverse    
sheaf on $Z_J$ is said to be a unipotent character sheaf (UCS) if it
appears in the perverse cohomology of $\p^w_!\bbq$ for some $w\in W$.
The UCS on $Z_J$ (up to isomorphism) form a set denoted by $UCS_J$.
We have $D=Z_I$ hence the set $UCS_I$ of UCS on $D$ is defined.

Let $J\sub I$ and let $(P,P',gU_P)\in Z_J$. Let
$z=\pos(P',P)\in{}^{\t(J)}W^J$. Following \cite{L04}, we define

(i) $(P^1,P'{}^1,gU_{P^1})$ by $P'{}^1=P'{}^P,P^1=g\i P'{}^Pg$.
\nl
(We use that $U_P\sub U_{P^1}$.) We have
$(P^1,P'{}^1,gU_{P^1})\in Z_{J_1}$ where\lb
$J_1=J\cap\t\i(zJz\i)$, see
\cite{L04,2.7(a)} and $\pos(P'{}^1,P^1)\in zW_J$, see \cite{L04, 3.2}.

If $z\in{}^{\t(J)}_*W$, we have $J_1=J$ and
$(P^1,P'{}^1,gU_{P^1})=(P,P',gU_P)$. If $z\in{}^{\t(J)}_\spa W$, we
have $\sha(J_1)<\sha(J)$.

\subhead 1.4\endsubhead
For $J\sub I$ we define a map $Z_J@>>>W$, $(P,P',gU_P)\m w_{P,P',gU_P}$,
by induction on $\sha(J)$.

Let $(P,P',gU_P)\in Z_J$. Let $z=\pos(P',P)\in{}^{\t(J)}W^J$.

If $z\in{}^{\t(J)}_*W$, we set $w_{P,P',gU_P}=z$. (In particular this
defines our map when $J=\em$.)

If $z\in{}^{\t(J)}_\spa W$, we set (with notation of 1.3(i))
$$w_{P,P',gU_P}=w_{P^1,P'{}^1,gU_{P^1}}.$$
(The right hand side is known by the induction hypothesis since
$(P^1,P'{}^1,gU_{P^1})\in Z_{J_1}$ with $\sha(J_1)<\sha(J)$).

By \cite{L06, 36.1}, for $(P,P',gU_P)\in Z_J$ we have

(a) $w=w_{P,P',gU_P}\in{}^{\t(J)}W$ and $\pos(P',P)=w_J$.

\subhead 1.5\endsubhead
For $J\sub I,w\in{}^{\t(J)}W$ we set 
$$Z_J^w=\{(P,P',gU_P)\in Z_J;w_{P,P',gU_P}=w\}.$$
This is a locally closed smooth subvariety of $Z_J$ stable under the
$G$-action $h:(P,P',gU_P)\m(hPh\i,hP'h\i,hgh\i U_{hPh\i})$. Note that
$(Z^w_J)_{w\in{}^{\t(J)}W}$ is a partition of $Z_J$.

If $J\sub I,w\in{}^JW$, then $(P,P',gU_P)\m(P^1,P'{}^1,gU_{P^1})$
(notation of 1.3(i)) defines a map $\vt_w:Z_J^w@>>>Z_{J_1}^w$
where $J_1=J\cap\t(w_JJw_J\i)$. This is an affine space bundle. See
\cite{L04,3.12}. If $w\in{}^{\t(J)}_*W$, this is the identity map.

If $w\in{}^{\t(J)}_*W$ we have (using 1.4(a)):

(a) $Z_J^w=\{(P,P',gU_P)\in Z_J;\pos(P',P)=w\}$.
\nl
In particular, $Z_J^w\ne\em$.

\subhead 1.6\endsubhead
Following \cite{L04}, for $J\sub I$, $w\in{}^{\t(J)}W$, we define the
notion of unipotent character sheaf (UCS) on $Z_J^w$ by induction on
$\sha(J)$.

Assume first that $w\in{}^{\t(J)}_*W$. Let $(P,P',gU_P)\in Z_J^w$. We
show that

(a) $P,P'$ have a common Levi subgroup.
\nl
(This was stated without proof in \cite{L04,4.6}.) Let
$(P^1,P'{}^1,gP_1)\in Z_{J_1}$ be as in 1.3(i), so that
$J_1=J\cap\t(wJw\i)=J$, see 1.1(a). Since $P'{}^1\sub P'$ and $J_1=J$,
it follows that $P'{}^1=P'$, that is $(P\cap P')U_{P'}=P'$. The obvious
homomorphism $P'@>>>P'/U_{P'}$ induces a surjective homomorphism\lb
$f:P\cap P'@>>>P'/U_{P'}$. Since $U_{P\cap P'}$ is a connected,
unipotent group, normal in $P\cap P'$, its image $f(U_{P\cap P'})$ is a
connected, unipotent group, normal in $P'/U_{P'}$. Since $P'/U_{P'}$ is
reductive, we see that $f(U_{P\cap P'})=\{1\}$. Thus $f$ induces a
surjective homomorphism $(P\cap P')/U_{P\cap P'}@>>>P'/U_{P'}$. Hence,
if $L$ is a Levi subgroup of $P\cap P'$, $f$ restricts to a surjective
homomorphism $L@>>>P'/U_{P'}$. It follows that $\dim(L)\ge\dim(L')$
where $L'$ is a Levi subgroup of $P'$ (which can be assumed to contain
$L$). This forces $L=L'$. Now $L$ is also contained in a Levi subgroup
$L''$ of $P$ which must have the same dimension as $L'$ hence as $L$.
It follows that $L=L''$. This proves (a).

Let $L$ be a common Levi subgroup of $P,P'$. Let
$D_w=\{g\in D;gPg\i=P',gLg\i=L\}$ (a connected component of the
reductive group \lb $NL=\{g\in\tG;gLg\i=L\}$ with identity component $L$.)
We consider the diagram
$$D_w@<j_1<<G\T D_w@>j_2>>Z_J^w,$$
where $j_1$ is the second projection and 
$$j_2(h,g)=(hPh\i,hP'h\i,hgh\i U_{hPh\i}).$$
Here $j_1$ is smooth with connected fibres and $j_2$ is an affine
space bundle. If $A$ is a UCS on $D_w$ (see 1.3), then
$j_1^\star(A)=j_2^\star(A')$ for a well defined simple perverse sheaf
$A'$ on $Z_J^w$. The set of such $A'$ (up to isomorphism) is denoted
by $UCS^w_J$. These are the UCS on $Z_J^w$.

We say that $A'\in UCS^w_J$ is cuspidal if the corresponding UCS on
$D_w$ is cuspidal. Let
$$UCS^w_{J,cu}=\{A'\in UCS^w_J;A'\text{ cuspidal}\}.$$

Assume next that $w\in{}^{\t(J)}_\spa W$. We set
$J_1=J\cap\t\i(w_JJw_J\i)$. We have $\sha(J_1)<\sha(J)$ and
$w\in{}^{\t(J_1)}W$ hence the UCS on $Z_{J_1}^w$ are defined.
Then the UCS on $Z_J^w$ are the simple perverse sheaves
$\vt_w^\star(A)$ (see 1.5) for various UCS $A$ on $Z_{J_1}^w$.
Their isomorphism classes form a set $UCS^w_J$.

\subhead 1.7\endsubhead
If $X$ is a variety over $\kk$ we denote by $\cd_m(X)$ the bounded
mixed derived category of $\bbq$-sheaves on $X$.
If $X$ has a given collection $UCS(X)$ of simple perverse sheaves 
we denote by $\cd_m^{UCS}(X)$ the subcategory of $\cd_m(X)$
consisting of complexes whose perverse cohomology sheaves have all
their simple subquotients in $UCS(X)$. For such $X$ we denote by
$\ck(X)$ the free $\ca$-module with basis given by $UCS(X)$. For
$C\in\cd_m^{UCS}(X)$ we set
$$gr(C)=\sum_{A\in UCS(X)}
\sum_{j,h\in\ZZ}(-1)^j(\text{multiplicity of } A \text{ in }
{}^pH^j(C)_h)v^hA\in\ck(X)$$
where the subscript $h$ denotes the subquotient of pure weight $h$
of a mixed perverse sheaf.

This applies in particular when $X$ is $Z_J$ or $Z_J^w$ with
$J\sub I$, $w\in{}^JW$ and hence also when $X=D$. (We have
$D=Z_I=Z_I^1$.)
We write $\ck_J$, $\ck^w_J$ instead of $\ck(X)$ when $X$ is $Z_J$ or
$Z_J^w$.
In each of these cases any UCS on $X$ is in $\cd_m^{UCS}(X)$ (it has
pure weight $0$).

Let $J\sub I$, $w\in{}^{\t(J)}W$. Let $i_{J,w}:Z_J^w@>>>Z_J$ be the
inclusion.
If $C\in\cd_m^{UCS}(Z_J)$ then $i_{J,w}^*C\in\cd_m^{UCS}(Z^w_J)$. See
\cite{L06, 36.3(a)}. If $C\in\cd_m^{UCS}(Z_J^w)$ then
$i_{J,w!}C\in\cd_m^{UCS}(Z_J)$. See \cite{L06, 36.3(f)}.
Hence we obtain an $\ca$-linear map $\ck_J@>>>\ck^w_J$ given by
$A\m gr(i_{J,w}^*A)$ for any $A\in UCS_J$ (we denote it again by
$i_{j,w}^*$) and an $\ca$-linear map $\ck^w_J@>>>\ck_J$ given by
$A\m gr(i_{J,w!}A)$ for any $A\in UCS^w_J$ (we denote it again by
$i_{j,w!}$).

By \cite{L06,36.8(a)} we have an $\ca$-linear isomorphism
$\ck_J@>>>\op_{w\in{}^{\t(J)}W}\ck_J^w$
given by $A\m\sum_w(i_{J,w}^*A)$ for any $A\in UCS_J$.
The inverse isomorphism restricted to $\ck_J^w$ is $A\m i_{J,w_!}A$
for any $A\in UCS^w_J$. We use this isomorphism to identify
$\ck_J=\op_{w\in{}^{\t(J)}W}\ck_J^w$. We have
$\ck_J=\ck_J^*\op\ck_J^\spa$ where

$\ck_J^*=\op_{w\in{}^{\t(J)}_*W}\ck_J^w$,
$\ck_J^\spa=\op_{w\in{}^{\t(J)}_\spa W}\ck_J^w$.

\subhead 1.8\endsubhead
Assume that $J'\sub J\sub I$. Let
$$\align&Z_{J',J}=\{(P,P',Q,Q',gU_Q);P\in\cp_{J'},P'\in\cp_{\t(J')},
Q\in\cp_J,\\&
Q'\in\cp_{\t(J)},P\sub Q,P'\sub Q',gU_Q\in D/U_Q,gPg\i=P'\}\endalign$$
(note that $U_Q\sub U_P$). We define $Z_{J'}@<c<<Z_{J',J}@>d>>Z_J$
by
$$c(P,P',Q,Q',gU_Q)=(P,P',gU_P), d(P,P',Q,Q',gU_Q)=(Q,Q',gU_Q).$$
Note that $c$ is an affine space bundle (the fibre at $(P,P',gU_P)$ is
isomorphic to $U_P/U_Q$ where $P\sub Q,P\in\cp_{J'},Q\in\cp_J$); $d$ is
proper (the fibre at $(Q,Q',gU_Q)$ is a partial flag manifold of $Q$).
Following \cite{L04,6.4}, we define
$$f_{J',J}=d_!c^*:\cd_m^{UCS}(Z_{J'})@>>>\cd_m^{UCS}(Z_J).$$
This induces an $\ca$-linear map $\ck_{J'}@>>>\ck_J$ given by
$A\m gr(f_{J',J}A)$ for any $A\in UCS_{J'}$ (we denote this map again
by $f_{J',J}$).

\subhead 1.9\endsubhead
For $J'\sub J$, $w\in{}^{\t(J)}W$, $y\in{}^{\t(J')}W$ we set
$$Z_{J',J}^{y,w}=\{(P,P',Q,Q',gU_Q)\in Z_{J',J};
(P,P',gU_P)\in Z_{J'}^y,(Q,Q',gU_Q)\in Z_J^w\}.$$

We define
$$Z_{J'}^y@<c_1<<Z_{J',J}^{y,w}@>d_1>>Z_J^w$$
where $c_1,d_1$ are the restrictions of $c,d$ in 1.8. Here $c_1$
is an affine space bundle.

Following \cite{L06} we define
$$f_{J',J}^{y,w}=d_{1!}c_1^*:\cd_m^{UCS}(Z_{J'}^y)@>>>
\cd_m^{UCS}(Z_J^w).$$
See \cite{L06, 36.5}. This induces an $\ca$-linear map
$\ck_{J'}^y@>>>\ck_J^w$ denoted again by $f_{J',J}^{y,w}$.

For $\x\in\ck^y_{J'}$ we have
$$f_{J',J}(\x)=\sum_{w\in{}^JW}f_{J',J}^{y,w}(\x).$$
See \cite{L06,36.5(a)}.

Let $J'\sub J$, $w\in{}^{\t(J)}_*W$, $y\in{}^{\t(J')}_*W$.

(a) If $Z_{J',J}^{y,w}\ne\em$ then $w=y_J$. See \cite{L06, 36.7(a)}.

\subhead 1.10\endsubhead
Assume that $K\sub J\sub J'\sub I$. Let $y\in{}^{\t(K)}_*W$,
$x\in{}^{\t(J'(}_*W$, $w\in{}^{\t(K)}_*W$ be such that
$x=y_{J'},w=y_J$. We show:

(a) $f_{K,J}^{y,w}=f_{J',J}^{x,w}f_{K,J'}^{y,x}$ as maps
$\cd_m(Z_K^y)@>>>\cd_m(Z_J^w)$.
\nl
We have the diagram
$$Z^y_K@<c<<Z_{K,J'}^{y,x}@>d>>Z^x_{J'}@<c'<<Z^{x,w}_{J',J}
@>d'>>Z^w_J$$
where $c,c'$ (resp. $d,d'$) are the analogues of $c_1$ (resp. $d_1$)
in 1.9.
The fibre product of $Z_{K,J'}^{y,x},Z_{J',J}^{x,w}$ over
$Z^x_{J'}$ (via $d,c'$) is
$$\align&Z_{K,J}^{y,w}=\{(S,S',Q,Q',gU_Q);S\in\cp_K,S'\in\cp_{\t(K)},
Q\in\cp_J,Q'\in\cp_{\t(J)},\\&S\sub Q, S'\sub Q',gU_Q\in D/U_Q,gSg\i=S',
\pos(S',S)=y\}.\endalign$$
We have a diagram
$$Z^y_K@<c<<Z_{K,J'}^{y,x}@<e<<Z_{K,J}^{y,w}
@>e'>>Z^{x,w}_{J',J}@>d'>>Z^w_J$$
where
$$e(S,S',Q,Q',gU_Q)=(S,S',P,P',gU_P),$$
$$e'(S,S',Q,Q',gU_Q)=(P,P',Q,Q',gU_Q).$$
(Here $P\in\cp_{J'}$ is defined by $S\sub P$ and $P'\in\cp_{\t(J')}$
is defined by $S'\sub P'$.) We have

(b) $d'_!c'{}^* d_!c^*=d'_!e'_!e^*c^*=(d'e')_!(ce)^*$
\nl
where $ce(S,S',Q,Q',gU_Q)=(S,S',gU_S)$,
$d'e'(S,S',Q,Q',gU_Q)=(Q,Q',gU_Q)$. Now (b) implies that (a) holds.

\subhead 1.11\endsubhead
In this subsection we assume that $\tG=D=G$ so that $\t=1$. Let
$J\sub I$. Let $NW_J=\{w\in W;wW_Jw\i=W_J\}$ be the normalizer of $W_J$
in $W$. Let $\cw_J$ be the set of elements of $NW_J$ that have minimal
length in their $W_J$-coset. We show:

(a) $\cw_J$ is a subgroup of $W$ and ${}^J_*W=\cw_J$.
\nl
By \cite{L76, 5.2(i)}, $\cw_J$ is a subgroup of $NW_J$  and
$\cw_J\sub{}^JW$. Moreover if $w\in\cw_J$ then clearly $w=w_J$ hence
$w\in{}^J_*W$. Thus $\cw_J\sub{}^J_*W$. Now let $w\in{}^J_*W$. Then
$w_JJw_J\i=J$. By 1.1(a) we have $w=w_J$ hence $w\in NW_J$. Since
$w\in{}^JW$ we see that $w\in \cw_J$. Thus ${}^J_*W=\cw_J$. This proves
(a).

\subhead 1.12\endsubhead
We still assume that $\tG=D=G$.

For $J'\sub I$, $J\sub I$, we write $J'\prec J$ whenever $J'\sub J$,
$J'\ne J$. For $J\sub I$ let
\define\tck{\ti{\ck}}
$\tck_J=\sum_{J'\prec J}f_{J',J}(\ck_{J'})$,
$\tck_J^*=\sum_{J'\prec J}f_{J',J}\ck^*_{J'}$.

By \cite{L06,36.6(d),(e)} we have

(a) $\ck^\spa_J\sub\tck_J$.
\nl
We show:

(b) $\tck_J=\tck^*_J$.
\nl
We argue by induction on $\sha(J)$. When $J=\em$, both sides of (b)
are $0$. Assume now that $J\ne\em$. Using (a) we have
$$\tck_J=\tck_J^*+\sum_{J'\prec J}f_{J',J}\ck^\spa_{J'}\sub
\tck^*_J+\sum_{J'\prec J}f_{J',J}\tck_{J'}.$$

By the induction hypothesis for $J'\prec J$ we have
$\tck_{J'}=\tck^*_{J'}$ hence

$$\align&\tck_J\sub\tck^*_J+\sum_{J'\prec J}f_{J',J}\tck^*_{J'}
=\tck^*_J+\sum_{J'\prec J,J''\prec J'}f_{J',J}f_{J'',J'}\ck^*_{J''}\\&
=\tck^*_J+\sum_{J'\prec J,J''\prec J'}f_{J'',J}\ck^*_{J''}=\tck^*_J.\endalign$$
Thus the left hand side of (b) is contained in the right hand side.
The opposite containment is obvious. This proves (b).

\head 2. Symplectic groups\endhead
\subhead 2.1\endsubhead
Let $C_n$ ($n\in\NN$) be the category whose objects are the
$\kk$-vector spaces $V$ of dimension $2n$ with a given nondegenerate symplectic form
$<,>:V\T V@>>>\kk$.
Let $\fC_n$ be the category whose objects are the objects $(V,<,>)$ of
$C_n$ with a given $F_q$-rational structure with Frobenius map
$F:V@>>>V$ such that $$<F(x),F(x')>=<x,x'>^q$$
for any $x,x'$ in $V$.
For $V\in C_n$ we denote by $Sp(V)$ the group of vector space
automorphisms of $V$ respecting $<,>$.
For $V\in\fC_n$ we define $F:Sp(V)@>>>Sp(V)$ by
$F(g)F(x)=F(gx)$ for $g\in Sp(V),x\in V$; then
$Sp(V)^F$ is a finite symplectic group.
If $V\in\fC_n$ with $n=k^2+k$ for some $k\in\NN$ we denote by
$\c(V):Sp(V)^F@>>>\ZZ$ the character of the unipotent cuspidal
representation of $Sp(V)^F$ and by $\ti\c(V):Sp(V)^F@>>>\QQ$ the almost
character associated to the cuspidal UCS of $Sp(V)$.
From the definition of almost characters we see that we have
$$\ti\c(V)=2^{-k}(\c(V)+\sum_{i=1}^{2^k-1}\d_i\c_i(V))\tag a$$
where for $i=1,2,\do,2^{2k}-1$, we have $\d_i\in\{1,-1\}$ and
$\c_i(V):Sp(V)^F@>>>\ZZ$ is a well
defined character of a unipotent noncuspidal irreducible representation
of $Sp(V)^F$.

\subhead 2.2\endsubhead
We now fix $V\in\fC_n$. For a subspace $A$ of $V$ we set
$A^\pe=\{x\in V;<x,A>=0\}$.

In this section we assume that $\tG=D=G=Sp(V)$, so that $\t=1$.
The Weyl group $W$ is of type $B_n$. (By convention, when $W=\{1\}$
we say that $W$ is of type $B_0$.) The simple reflections of $W$ form
a set $I=\{s_1,s_2,\do,s_n\}$ in which $s_is_{i+1}$ has order $3$ for
$i=1,2,\do,n-2$ and order $4$ for $i=n-1$.

For $t\in[0,n]$ let $E_t$ be the set of all isotropic subspaces of
dimension $t$ of $V$. Let $(A,A')\in E_t\T E_t$. We set
$$[A,A']=(A^\pe\cap A'{}^\pe)/(A\cap A').$$
We show:

(a)  $\dim[A,A']=2n-2t$.
\nl
Indeed, $[A,A']$ is $(A+A')^\pe/(A\cap A')$ hence it has dimension
$2n-\dim(A+A')-\dim(A\cap A')=2n-\dim(A)-\dim(A')=2n-2t$.

We say that $(A,A')$ is {\it good } if there exists a subspace $L$
of $V$ such that $A^\pe=L\op A,A'{}^\pe=L\op A'$. We show:

(b) the following three conditions are equivalent:

(i) $(A,A')$ is good;

(ii) $A\cap A'=A^\pe\cap A'$;

(iii) $A\cap A'=A\cap A'{}^\pe$.
\nl
Assume that (i) holds. We show that (ii) holds. Let $L$ be a subspace
of $V$ such that $A^\pe=L\op A$, $A'{}^\pe=L\op A'$. Let
$a'\in A^\pe\cap A'$.
We have $l=a'-a$ where $l\in L,a\in A,a'\in A'$. For any $l'\in L$
we have $<l,l'>=<a',l'>-<a,l'>=0$ (since $l'\in A^\pe\cap A'{}^\pe$).
Since $<,>$ is nondegenerate on $L$, it follows that $l=0$ so that
$a'=a\in A\cap A'$. Thus, $A^\pe\cap A'\sub A\cap A'$; the reverse
inclusion is obvious hence (ii) holds. Similarly (iii) holds (since
(i) is symmetric in $A,A'$).

Assume that (ii) holds. We show that (i) holds. Let $L$ be a subspace
of $V$ such that $A^\pe\cap A'{}^\pe=L\op(A\cap A')$. The obvious
linear map $\g_{A,A'}:[A,A']@>>>A^\pe/A$ is injective (by (ii)).
Since $\dim(A^\pe/A)=\dim(L)=2n-2t$ (see (a)), this map is an
isomorphism of vector spaces compatible with the symplectic forms
induced by $<,>$; it follows that $A^\pe=L\op A$. In particular, the
symplectic form on $[A,A']$ is nondegenerate. The obvious linear map
$[A,A']@>>>A'{}^\pe/A'$ is a linear map compatible with the symplectic
forms; since these forms are nondegenerate and the vector spaces
involved have the same dimension, this map must be an isomorphism; it
follows that $A'{}^\pe=L\op A'$. Thus, (i) holds.

Assume that (iii) holds. Then (i) holds. Indeed this follows from the
previous argument by interchanging $A,A'$ (note that (i) is symmetric
in $A,A'$). This proves (b).

One can show that conditions (i)-(iii) are equivalent to

(iv) $A\cap A'=A^\pe\cap A'{}^\pe\cap(A+A')$.
\nl
From the proof of (b) we see that when $(A,A')$ is good,

(c) the obvious linear map $\g_{A,A'}:[A,A']@>>>A^\pe/A$ is an
isomorphism (in $C_{n-t}$).
\nl
Note that if $t$ is $0$ or $n$, then any $(A,A')\in E_t\T E_t$ is good.
If $t=1<n$, then $(A,A')\in E_t\T E_t$ is good if $A=A'$ or if
$A'\cap A^\pe=0$ and is not good if $A\ne A'$ and $A'\sub A^\pe$.

\subhead 2.3\endsubhead
Assume that $(A,A')$ is good.
There is a unique isomorphism $\psi^A_{A'}:A^\pe/A@>>>A'{}^\pe/A'$
(in $C_{n-t}$) such that the composition 
$$[A,A']@>\g_{A,A'}>>A^\pe/A@>\psi^A_{A'}>>A'{}^\pe/A'$$
is equal to $[A',A]@>\g_{A',A}>>A'{}^\pe/A'$.
(We use that $[A,A']=[A',A]$; $\g_{A,A'}$ is the isomorphism
in 2.2(c).) We have $\psi^A_{A'}\psi^{A'}_A=1$.

We define $F:E_t@>>>E_t$ by $A\m F(A)$; then
$E_t^F=\{A\in E_t;F(A)=A\}$ is defined.
If $A\in E_t^F$ then $A^\pe/A$ is naturally an object of $\fC_{n-t}$.
If $(A,A')\in E_t^F\T E_t^F$ is good then $\psi^A_{A'}$ is an
isomorphism in $\fC_{n-t}$.

\subhead 2.4\endsubhead
Let $\ce_t$ be the set of all sequences
$V_*=(V_1\sub V_2\sub\do\sub V_t)$ with $V_i\in E_i$ for $i=1,2,\do,t$.
For $V_*$ as above let
$$P(V_*)=\{g\in G;gV_i=V_i\text{ for }i=1,2,\do,t\}.$$
This is a parabolic subgroup of $G$ in $\cp_{J_t}$ where
$J_t=\{s_{t+1},\do,s_{n-1},s_n\}$. Its unipotent radical is denoted by
$U(V_*)$.
We can identify $Z_{J_t}$ with the set of all triples
$(V_*,V'_*,gU(V_*))$ where $V_*\in\ce_t,V'_*\in\ce_t$, $g\in G$,
$g(V_*)=V'_*$, by writing $(V_*,V'_*,gU(V_*))$ instead of
$(P(V_*),P(V'_*),gU(V_*))$.

Let $(V_*,V'_*,gU(V_*))\in Z_{J_t}$ be such that $(V_t,V'_t)$ is
good (or equivalently such that $P(V_*),P(V'_*)$ contain a common
Levi subgroup.)
Let $\un g:V_t^\pe/V_t@>>>V'_t{}^\pe/V'_t$ be the isomorphism
induced by $g:V@>>>V$ (which carries $V_t$ onto $V'_t$ and
$V_t^\pe$ onto $V'_t{}^\pe$). We set
$$\bg=\psi^{V'_t}_{V_t}\un g:V_t^\pe/V_t@>>>V_t^\pe/V_t.$$
(Notation of 2.3.) This is an element of $Sp(V_t^\pe/V_t)$ which
depends only on $gU(V_*)$.

We define $F:\ce_t@>>>\ce_t$ by
$(V_1\sub V_2\sub\do\sub V_t)\m(F(V_1)\sub F(V_2)\sub\do\sub F(V_t))$;
then $\ce_t^F$ is defined.
We define $F:Z_{J_t}@>>>Z_{J_t}$ by
$$(V_*,V'_*,gU(V_*))\m(F(V_*),F(V'_*),F(g)U(F(V_*)))$$
so that
$$(Z_{J_t})^F=\{(V_*,V'_*,gU(V_*))\in Z_{J_t};V_*\in\ce_t^F,
V'_*\in\ce_t^F,g\in Sp(V)^F)\}$$
is defined.

\subhead 2.5\endsubhead
If $t\in[0,n]$ then $W_{J_t}$ is a Weyl group of type $B_{n-t}$ and
$\cw_{J_t}$ is a Weyl group of type $B_t$ with generators
$s_1,s_2,\do,s_{t-1},s'_t$ where $s'_t$ is the longest
element of $W$ times the longest element of $W_{J_t}$. 

If $t\in[0,n],t'\in[0,n],t+t'\in[0,n]$, we define
$\cw_{J_{t,t'}}$ to be the subgroup of $W_t$ with generators
$s_{t+1},s_{t+2},\do,s_{t+t'-1},s'_{t+t'}$ where $s'_{t+t'}$ is the
longest element of $W_{J_t}$ times the longest element of
$W_{J_{t+t'}}$; this is a Weyl group of type $B_{t'}$ commuting with
$\cw_{J_t}$. Note that the product $\cw_{J_t}\cw_{J_{t,t'}}$ is
contained in $\cw_{J_{t+t'}}$.)

Let $\Xi$ be the set of all $t\in[0,n]$ such that $n-t=k^2+k$ for
some $k\in\NN$. 

We set
$$M=\{(w,t);t\in\Xi,w\in\cw_{J_t}\}.$$

For $(w,t)\in M$ let 
$$\YY^w_t=\{(V_*,V'_*,g)\in\ce_t\T\ce_t\T G;g(V_*)=V'_*,
\pos(V'_*,V_*)=w\}.$$
(We write $\pos(V'_*,V_*)$ instead of $\pos(P(V'_*),P(V_*))$.)
Note that $\YY^w_t$ has an obvious Frobenius map
$F:\YY^w_t@>>>\YY^w_t$ with fixed point set $(\YY^w_t)^F$.

Let $t\in\Xi$. Let $\Irr\cw_{J_t}$ be the set of (isomorphism classes
of) irreducible representations (over $\QQ$) of $\cw_{J_t}$.

Let $V_*\in\ce_t^F$. Then $Sp(V_t^\pe/V_t)^F)$ is naturaly a quotient of
$P(V_*)^F:=P(V_*)\cap G^F$. Hence the unipotent cuspidal
representation of $Sp(V_t^\pe/V_t)^F)$ (over $\QQ$) can be regarded as a
representation of $P(V_*)^F$ (over $\QQ$). We induce this representation
from $P(V_*)^F$ to $G^F$. We obtain a representation
$\VV_t$ (over $\QQ$) of $G^F$. From \cite{L77} it is known that
the algebra $\End_{G^F}\VV_t$ is an Iwahori-Hecke algebra $H_q$ (over
$\QQ$) with possibly unequal parameters associated to the Weyl group
$\cw_{J_t}$; this algebra has a $\QQ$-basis $\{T_w;w\in\cw_{J_t}\}$
described in \it{loc.cit.} Moreover, as a $H_q\T G^F$-module we have
$$\VV_t=\op_E E_q\ot\tE_q$$ where $E$ runs over $\Irr\cw_{J_t}$ and for
each such $E$, $E_q$ denotes the irreducible representation of $H_q$
(over $\QQ$) corresponding to $E$, and $\tE_q$ denotes the irreducible
representation of $G^F$ (over $\QQ$) corresponding to $E_q$.

Let $\cv_t$ be the $\QQ$-vector space 
 with basis given by the functions $G^F@>>>\QQ$,
$g\m\tr(g,\tE_q)$ (for various $E\in\Irr\cw_{J_t}$).
We have
$$\dim\cv_t=\sha(\Irr\cw_{J_t})=\sha(\cw^!_{J_t})$$
where $\cw^!_{J_t}$ is a set of representatives for the conjugacy
classes in $\cw_{J_t}$.

For $w\in\cw_{J_t}$ we define $Y^w_t:G^F@>>>\ZZ$ by
$$Y^w_t(g)=\sum_{(V_*,V'_*,g)\in(\YY^w_t)^F}\c(V_t^\pe/V_t)(\bg).$$
(Here $\c(V_t^\pe/V_t)$ is as in 2.1 and $\bg$ is as in 2.4).

From the definition we see that for $w\in\cw_{J_t},g\in G^F$, we have
$$Y^w_t(g)=\tr(T_wg,\VV_t)$$
hence
$$Y^w_t(g)=\sum_{E\in\Irr\cw_{J_t}}\tr(T_w,E_q)\tr(g,\tE_q).$$
Note that $Y^w_t\in\cv_t$.

From the previous equality we can deduce that for any
$E\in\Irr\cw_{J_t}$ we have
$$\tr(g,\tE_q)=\sum_{w\in\cw^!_{J_t}}m_{w,E}(q)Y^w_t(g)\tag a$$
where $m_{w,E}(q)$ is the value at $q=v^2$ of an element
$m_{w,E}\in\QQ(v^2)$ (at least if $q$ is sufficiently large). Indeed,
the square matrix with entries $\tr(T_w,E_q)$,
$w\in\cw^!_{J_t},E\in\Irr\cw_{J_t}$,
can be regarded as a matrix with entries on $\ca$ (specialized with
$v=q^{1/2}$) whose determinant is in $\ca-\{0\}$ hence this
determinant at $v=q^{1/2}$ is $\ne0$ if $q$ is sufficiently large.
(The specialization for $v=1$ of the last matrix is the matrix with
entries $\tr(w,E)$, $w\in\cw^!_{J_t},E\in\Irr\cw_{J_t}$, which is well
known to be invertible.)

{\it In the remainder of this paper we assume that $q$ is sufficiently
large.}

From (a) we see that

(b) $\{Y^w_t;w\in\cw^!_{J_t}\}$ is a $\QQ$-basis of $\cv_t$.
\nl
Let $\cv$ be the $\QQ$-vector space of functions $G^F@>>>\QQ$ with
basis given by the characters of the unipotent representations of $G^F$.
By \cite{L77} we have $\cv=\op_{t\in\Xi}\cv_t$.
Note that $Y^w_t\in\cv$ for any $(w,t)\in M$); from (b) we see that

(c) $\{Y^w_t;t\in\Xi,w\in\cw^!_{J_t}\}$ is a $\QQ$-basis of $\cv$.

\subhead 2.6\endsubhead
Let $t\in\Xi$ (so that $t=n-(k^2+k)$ with $k\in\NN$).
For $w\in\cw_{J_t}$ we define ${}'Y^w_t:G^F@>>>\QQ$ by
$${}'Y^w_t(g)=\sum_{(V_*,V'_*,g)\in(\YY^w_t)^F}\ti\c(V_t^\pe/V_t)(\bg).$$
(Here $\ti\c(V_t^\pe/V_t)$ is as in 2.1 and $\bg$ is as in 2.4.)
Using 2.1(a) we see that
$${}'Y^w_t(g)=2^{-k}(Y^w_t(g)+\sum_{i=1}^{2^k-1}\d_i{}^iY^w_t(g))
\tag a$$
where $\d_i\in\{1,-1\}$ and
$${}^iY^w_t(g)=\sum_{(V_*,V'_*,g)\in(\YY^w_t)^F}\c_i(V_t^\pe/V_t)(\bg).
$$
Using 2.5(a) applied to $V_t^\pe/V_t$ instead of $V$ we have
$$\c_i(V_t^\pe/V_t)=\sum_{y\in\cw_{J_{t,t'}}^!}m_{y,i}(q)Y^y_{t'}$$
where $t'=k^2+k-(k'{}^2+k')>0$, $\cw_{J_{t,t'}}^!$ is a set of
representatives for the conjugacy classes in $\cw_{J_{t,t'}}$ (see 2.5)
and $m_{y,i}(q)$ is the value at $q=v^2$ of an element
$m_{y,i}\in\QQ(v^2)$. Here $Y^y_{t'}$ is the analogue
for $V_t^\pe/V_t$ instead of $V$ of $Y^w_t$ in 2.5. It follows that
$${}^iY^w_t(g)=\sum_{(V_*,V'_*,g)\in(\YY^w_t)^F}
\sum_{y\in\cw_{J_{t,t'}}^!}m_{y,i}(q)Y^{wy}_{t+t'}(g).\tag b$$

Let $M^!=\{(w,t)\in M;w\in\cw^!_{J_t}\}$.
From (a),(b) we deduce for $w\in\cw_{J_t}$:
$${}'Y^w_t=2^{-k}Y^w_t+
\sum_{(w'',t'')\in M^!;t''>t}\mu_{w'',t''}(q)Y^{w''}_{t''}\tag c$$
where $\mu_{w'',t''}(q)$ is the value at $q=v^2$ of an element
$\mu_{w'',t''}\in\QQ(v^2)$. In particular,

(d) ${}'Y^w_t\in\cv$ for any $(w,t)\in M$.
\nl
By (c), the elements ${}'Y^w_t,(w,t)\in M$ are related to the basis
$Y^w_t,(w,t)\in M$ by an upper triangular matrix with nonzero
diagonal entries. It follows that

(e) $\{{}'Y^w_t;(w,t)\in M^!\}$ is a $\QQ$-basis of $\cv$ and
for any $(w,t)\in M^!$ we have
$$Y^w_t=2^k{}'Y^w_t+\sum_{(w'',t'')\in M^!;t''>t}\ti\mu_{w'',t''}(q)
{}'Y^{w''}_{t''}$$

where $\ti\mu_{w'',t''}(q)$ is the value at $q=v^2$ of an element
$\ti\mu_{w'',t''}\in\QQ(v^2)$. 

Let $\cv_{>0}=\op_{t\in\Xi;t>0}\cv_t$. We show:

(f) ${}'Y^w_t\in\cv_{>0}$ for any $(w,t)\in M-\{1,0\}$.
\nl
This follows from (c) since $Y^w_t\in\cv_{>0}$ for any
$(w,t)\in M-\{1,0\}$.

Using (e) and (f) we deduce

(g) $\{{}'Y^w_t;(w,t)\in M^!-\{1,0\}\}$ is a $\QQ$-basis of $\cv_{>0}$.

\subhead 2.7\endsubhead
Let $t\in[0,n]$. We fix $w\in\cw_{J_t}$.
Let $\Xi_t$ be the set of all $t'\in[0,n-t]$ such that
$t+t'\in\Xi$ that is $n-t-t'=k'{}^2+k'$ for some $k'\in\NN$. 
Let $M_t$ be the set of all pairs $(y,t')$ where $t'\in\Xi_t$ and
$y\in\cw_{J_{t,t'}}$ (see 2.5).
Let $M^!_t$ be the set of all $(y,t')$ in $M_t$ such that
$y\in\cw^!_{J_{t,t'}}$.

For $V_*=(V_1\sub V_2\sub\do\sub V_{t+t'})\in\ce_{t+t'}$ we set
$V_{*,t}=(V_1\sub V_2\sub\do\sub V_t)\in\ce_t$.

For $(y,t')\in M_t$ we set
$$\align&\YY^{w,y}_{t,t'}\\&=\{(V_*,V'_*,gU_{V_{*,t}});V_*\in\ce_{t+t'},
V'_*\in\ce_{t+t'},g\in G,g(V_*)=V'_*,\pos(V'_*,V_*)=wy\}.\endalign$$
Note that $\YY^{w,y}_{t,t'}$ has an obvious Frobenius map
$F:\YY^{w,y}_{t,t'}@>>>\YY^{w,y}_{t,t'}$ with fixed point set
$(\YY^{w,y}_{t,t'})^F$.

For $(y,t')\in M_t$ we define $Y^{w,y}_{t,t'}:(Z^w_{J_t})^F@>>>\ZZ$ by
$$Y^{w,y}_{t,t'}(\tV_*,\tV'_*,gU(\tV_*))=
\sum_{(V_*,V'_*,gU_{V_{*,t}})\in(\YY^{w,y}_{t,t'})^F;V_{*,t}=\tV_*,
V'_{*,t}=\tV'_*}\c(V_{t+t'}^\pe/V_{t+t'})(\bg)$$
and ${}'Y^{w,y}_{t,t'}:(Z^w_{J_t})^F@>>>\ZZ$ by
$${}'Y^{w,y}_{t,t'}(\tV_*,\tV'_*,gU(\tV_*))=
\sum_{(V_*,V'_*,gU_{V_{*,t}})\in(\YY^{w,y}_{t,t'})^F;V_{*,t}=\tV_*,
V'_{*,t}=\tV'_*}\tc(V_{t+t'}^\pe/V_{t+t'})(\bg).$$

We choose $$({}^0V_*,{}^0V'_*)\in\ce_t^F\T\ce_t^F$$ such that
$\pos({}^0V'_*,{}^0V_*)=w$.
Let $\ph:Sp(({}^0V_t^\pe/{}^0V_t)^F)@>>>\QQ$ be a class function.
There is a well defined function $\ti\ph:(Z^w_{J_t})^F@>>>\QQ$
such that $\ti\ph$ is constant on the orbits of the $G^F$-action
$$h:(\tV_*,\tV'_*,gU(\tV_*))\m(h(\tV_*(,h(\tV'_*),hgh\i U(h(\tV_*)))$$
on $(Z^w_{J_t})^F$ and such that
$\ti\ph(\tV_*,\tV'_*,gU(\tV_*))=\ph(\bg)$
if $\tV_*={}^0V_*,\tV'_*={}^0V'_*$ (with $\bg$ as in 2.4).
Now $\ph\m\ti\ph$ defines an injective linear map from the vector space of class functions $Sp({}^0V_t^\pe/{}^0V_t)^F@>>>\QQ$ to the vector
space of functions $(Z^w_{J_t})^F@>>>\QQ$.

From the definitions we see that $Y^{w,y}_{t,t'}$ (resp.
${}'Y^{w,y}_{t,t'}$) is the image under $\ph\m\ti\ph$ of the class
function $Y^y_{t'}$ (resp. ${}'Y^y_{t'}$) on 
$Sp({}^0V_t^\pe/{}^0V_t)^F$.
Hence from 2.6 we can deduce the following results.

(a) The vector space spanned by the functions $Y^{w,y}_{t,t'}$ (with
$(y,t')\in M_t$) on $(Z^w_{J_t})^F$ coincides with the vector space
spanned by the functions ${}'Y^{w,y}_{t,t'}$ (with $(y,t')\in M_t$) on
$(Z^w_{J_t})^F$.

(b) The functions $Y^{w,y}_{t,t'}$ (with
$(y,t')\in M^!_t$) form a basis of the vector space in (a). The
functions ${}'Y^{w,y}_{t,t'}$ (with $(y,t')\in M^!_t$) also form a
basis of the vector space in (a).

(c) We have     
$$Y^{w,1}_{t,0}=2^{k'}{}'Y^{w,1}_{t,0}
+\sum_{(y,t')\in M_t^!-\{1,0\}}\ti\mu_{y,t'}(q){}'Y^{w,y}_{t,t'}$$
where $\ti\mu_{y,t'}(v^2)\in\QQ(v^2)$ is as in 2.6(e).

(d) The $\QQ$-vector space spanned by the functions
$Y^{w,y}_{t,t'}$ (with $(y,t')\in M_t-\{0,1\}$) has a basis
consisting of $Y^{w,y}_{t,t'}$ (with $(y,t')\in M^!_t-\{0,1\}$).

\subhead 2.8\endsubhead
In the setup of 2.7 let $(y,t')\in M_t$. Let
$$\bar\YY^{w,y}_{t,t'}=\{(V_*,V'_*)\in\ce_{t+t'}\T\ce_{t+t'};
\pos(V'_*,V_*)=wy\}.$$
Let $\p:\YY^{w,y}_{t,t'}@>>>\bar\YY^{w,y}_{t,t'}$ be the obvious map.
Now $G$ acts on $\YY^{w,y}_{t,t'}$ by
$$h:(V_*,V'_*,gU_{V_1\sub\do\sub V_t})\m
(h(V_*),h(V'_*),hgh\i U_{h(V_1)\sub\do\sub h(V_t)}).$$
There is a (unique up to isomorphism) $G$-equivariant constructible
$\bbq$-sheaf $\cl$ on $\YY^{w,y}_{t,t'}$ such that for any
$(V_*,V'_*)\in\bar\YY^{w,y}_{t,t'}$, the restriction of $\cl$ to
$\p\i(V_*,V'_*)$, identified with
$$\{gU_{V_{*,t}},g\in G;g(V_*)=V'_*\},$$
is the inverse image of the constructible sheaf associated to a
cuspidal UCS on $Sp(V_{t+t'}^\pe/V_{t+t'})$ under the map
$g\m\bg$ (as in 2.4).

Let $\cl^{w,y}_{t,t'}=a_!\cl$ where $a:\YY^{w,y}_{t,t'}@>>>Z^w_{J_t}$
is
$$(V_*,V'_*,gU_{V_{*,t}})\m(V_{*,t},V'_{*,t},gU_{V_{*,t}}).$$
We can find an isomorphism $\a:F^*\cl@>>>\cl$ (unique up to a nonzero
scalar). From $\a$ we get an isomorphism
$F^*\cl^{w,y}_{t,t'}@>>>\cl^{w,y}_{t,t'}$.
Taking alternating sum of its traces on stalks of cohomology sheaves
at elements of $(Z^w_{J_t})^F$ we obtain the ``characteristic
function'' $\x_{\cl^{w,y}_{t,t'},F}:(Z^w_{J_t})^F@>>>\bbq$.

Since $q$ is sufficiently large, the results of \cite{L92} (connecting
cuspidal UCS with almost characters) are applicable. Using \cite{L92}
and the definitions we see that
$\x_{\cl^{w,y}_{t,t'},F}=\l{}'Y^{w,y}_{t,t'}$ for some
$\l\in\bbq^*$. By choosing $\a$ appropriately, we can achieve that
$\l=1$ so that $\x_{\cl^{w,y}_{t,t'},F}={}'Y^{w,y}_{t,t';q}$.
(We now write ${}'Y^{w,w'}_{t,t';q}$ instead of ${}'Y^{w,y}_{t,t'}$
to emphasize that it depends on $F_q$.)
Replacing $F$ by a power $F^s,s\ge1$ we have similarly
$$\x_{\cl^{w,y}_{t,t'},F^s}={}'Y^{w,y}_{t,t';q^s}:
(Z^w_{J_t})^{F^s}@>>>\QQ.$$

From the definitions we have
$\cl^{w,y}_{t,t'}\in\cd_m^{UCS}(Z^w_{J_t})$ and

(a) $gr(\cl^{w,y}_{t,t'})$ is (up to a factor $\pm v^r$) the image
under $f_{J_{t+t'},J_t}^{wy,w}$ of the unique object in
$UCS^{wy}_{J_{t+t'},cu}$.

\subhead 2.9\endsubhead
Let $(w,t)\in M$. We show:

(a) The vectors $\{gr(\cl^{w,y}_{t,t'});(y,t')\in M_t^!\}$ in
$\un\ck(Z^w_{J_t})$ are linearly independent over $\QQ(v)$.
\nl
Assume that this is not so. Then we can find $\phi_{y,t'}(v)\in\ca$
for $(y,t')\in M_t^!$, not all zero, such that
$$\sum_{(y,t')\in M^!_t}\phi_{y,t'}(v)gr(\cl^{w,y}_{t,t'})=0$$
in $\un\ck(Z_{J_t}^w)$. From this we deduce
$\sum_{(y,t')\in M^!_t}\phi_{y,t'}(q^s){}'Y^{w,y}_{t,t',q^s}=0$ 
for any $s\ge1$.
Since $\{{}'Y^{w,y}_{t,t',q^s};(y,t')\in M^!_t\}$ are linearly
independent for any $s\ge1$ (see 2.7), we see that $\phi_{y,t'}(q^s)=0$ for any
$(y,t')\in M^!_t$ and any $s\ge1$. It follows that $\phi_{y,t'}(v)=0$
for any $(y,t')\in M^!_t$. This contradiction proves (a).

Next we note that from the classification of UCS on $Z_{J_t}^w$ (see
\cite{L04}) it is known that
$$\dim_{\QQ(v)}\un\ck(Z_{J_t}^w)=\sha(M^!_t).$$
Using (a) it follows that

(b) the vectors $\{gr(\cl^{w,y}_{t,t'});(y,t')\in M^!_t\}$ form a
$\QQ(v)$-basis of $\un\ck(Z_{J_t}^w)$. In particular,
$\un\ck(Z_{J_t}^w)$ is spanned by $gr(\cl^{w,y}_{t,t'})$ for various
$(w,t)\in M$.
\nl
Let $\un\ck(Z_{J_t}^w)_0$ be the $\QQ(v)$-subspace of
$\un\ck(Z_{J_t}^w)$ spanned by $gr(\cl^{w,1}_{t,0}$.
Let $\un\ck(Z_{J_t}^w)_{>0}$ be the $\QQ(v)$-subspace of
$\un\ck(Z_{J_t}^w)$ spanned by $gr(\cl^{w,y}_{t,t'}$ for various
$(y,t')\in M_t-\{1,0\}$. By (b) we have
$$\un\ck(Z_{J_t}^w)=\un\ck(Z_{J_t}^w)_0+\un\ck(Z_{J_t}^w)_{>0}.$$ 
We show:

(c) $\un\ck(Z_{J_t}^w)=\un\ck(Z_{J_t}^w)_0\op\un\ck(Z_{J_t}^w)_{>0}$.
\nl
Assume that (c) does not hold. Then
$$gr(\cl^{w,1}_{t,0})\in\un\ck(Z_{J_t}^w)_{>0}.$$
It follows that 
$$\phi gr(\cl^{w,1}_{t,0})=\sum_{(y,t')\in M_t,t'>0}
\phi_{y,t'}gr(\cl^{w,y}_{t,t'})$$
where $\phi\in\ca-\{0\}$ and $\phi_{y,t'}\in\ca$ for
$(y,t')\in M_t-\{1,0\}$.
From this we deduce that for any $s\ge1$,
$\phi(q^s){}'Y^{w,1}_{t,0,q^s}$ is a linear combination of functions
${}'Y^{w,y}_{t,t',q^s}$ with $(y,t')\in M_t-\{1,0\}$ hence, using
2.7(d), it is a linear combination of functions
${}'Y^{w,y}_{t,t',q^s}$ with $(y,t')\in M^!_t-\{1,0\}$.
Using now 2.7(b) we deduce that $\phi(q^s)=0$.
Since this holds for any $s\ge1$ it follows that $\phi=0$.
This contradiction proves (c).

Let $K\sub I,y\in{}^K_*W$; let $\ck^y_{K,cu}$ be the
$\ca$-submodule of $\ck^y_K$ with basis $UCS^y_{K,cu}$  (see 1.6).

Using 2.9(c) and 2.8(a) we see that the following holds.

\proclaim{Proposition 2.10}Let $t\in[0,n]$. Let $w\in\cw_{J_t}$. We have
$$\un\ck^w_{J_t}=\un\ck^w_{J_t,cu}\op\sum_{(y,t')\in M_t-\{1,0\}}
f_{J_{t+t'},J_t}^{wy,w}\un\ck^{wy}_{J_{t+t'},cu}.$$
\endproclaim

\subhead 2.11\endsubhead
Let $J\sub I$. Let $t$ be the smallest number in $[0,n]$ such that
$J_t\sub J$. We show:

(a) The linear map
$$\op_{y\in\cw_{J_t}}\un\ck_{J_t}^y@>>>\op_{w\in\cw_J}\un\ck_J^w$$
with components $f_{J_t,J}^{y,w}$ is surjective.
\nl
Let $L'$ be a Levi subgroup of some $P'\in\cp_{J_t}$; define
$P''\in\cp_J$ by $P'\sub P''$ and let $L''$ be the Levi subgroup of
$P''$ such that $L'\sub L''$. From the definitions, the map in (a) can
be identified with a linear map
$$\op_{D'}\un\ck(D')@>>>\op_{D''}\un\ck(D'')\tag b$$
where $D'$ (resp. $D''$) runs over the set of connected components of
the normalizer $NL'$ of $L'$ (resp. $NL''$ of $L''$) in $G$. Let $\ct$
be the connected centre of $L'$. We have $L'/\ct=L'_0$ where
$L'_0=Sp(Y)$ with $Y$ a subspace of $V$ of dimension $2(n-t)$ on which
$<,>$ is nondegenerate; we have canonically $NL'/\ct=L'_0\T W(Y^\pe)$
where $W(Y^\pe)$ is the Weyl group of $Sp(Y^\pe)$.
Now $\ct$ is also central in $L''$ and we have $L''/\ct=L'_0\T L''_0$
where $L''_0$ is a product of general linear groups contained in
$Sp(Y^\pe)$; we have canonically $NL''/\ct=L'_0\T N_{Sp(Y^\pe)}(L''_0)$.
We have canonically $UCS(D')=UCS(D'/\ct)$, $UCS(D'')=UCS(D''/\ct)$ for
any $D',D''$ in (b). Hence (b) can be rewritten as a linear map
$$\op_{z\in W(Y^\pe)}\un\ck(L'_0)\ot\bbq^z@>>>
\op_\d\un\ck(L'_0)\ot\ck(\d)$$

where $\bbq^z$ is a copy of $\bbq$ corresponding to $z$ and $\d$ runs
through the set of connected components of  $N_{Sp(Y^\pe)}(L''_0)$.
This map is obtained by tensoring the identity map of $\un\ck(L'_0)$
with a linear map
$$\op_{z\in W(Y^\pe)}\bbq^z@>>>\op_\d\ck(\d).$$

Replacing $Y^\pe$ by $V$, we are reduced to the case where
$t=n$ so that $J_t=\em$ and $W_J$ is a product of Weyl groups of type
$A$. The last linear map becomes a linear map
$$\op_{z\in W}\bbq^z@>>>\op_{w\in\cw_J}\ck(\d_w)=\op_{w\in\cw_J}\ck_J^w
$$

where $\d_w$ is the connected component of $NL''$ corresponding to $w$
under $NL''/l=L''\lra\cw_J$ ($L''$ is now a Levi subgroup of some
$P''\in\cp_J$). It is enough to show that this linear map is surjective.
This follows from the fact that for any $w\in\cw_J$, the linear map 
$$\op_{z\in W;z\in W_Jw=wW_J}\bbq^z@>>>\ck(\d_w)=\ck_J^w$$

with components $f_{\em,J}^{z,w}$ is surjective; this is a known
property of UCS on connected components of (possibly disconnected)
groups whose identity component is a product of general linear groups.

Returning to the setup of (a) we have the following consequence of (a)
and 1.9(a).

(c) For any $w\in\cw_J$, the linear map
$$\op_{y\in\cw_{J_t};y_J=w}\un\ck_{J_t}^y@>>>\un\ck_J^w$$
with components $f_{J_t,J}^{y,w}$ is surjective.

\proclaim{Proposition 2.12} Let $J\sub I$. Let $w\in\cw_J$. We have
$$\un\ck^w_J=\sum_{(z,s)\in M;J_s\sub J;w=z_J}f_{J_s,J}^{z,w}
\un\ck^z_{J_s,cu}.$$
\endproclaim
Define $t\in[0,n]$ as in 2.11. By 2.11(a) we have
$$\un\ck_J^w=\sum_{y\in\cw_{J_t};y_J=w}f_{J_t,J}^{y,w}\un\ck_{J_t}^y.$$
By 2.10, for any $y$ in the sum we have
$$\un\ck_{J_t}^y\sub\sum_{s\in\Xi,s\ge t,z\in\cw_{J_s}}f_{J_s,J_t}^{z,y}
\ck^z_{J_s,cu}$$
hence
$$\un\ck_J^w\sub\sum_{s\in\Xi,s\ge t,z\in\cw_{J_s},y\in\cw_{J_t};y_J=w}
f_{J_t,J}^{y,w}f_{J_s,J_t}^{z,y}\ck^z_{J_s,cu}$$
Using 1.10(a) we obtain
$$\un\ck_J^w\sub\sum_{(z,s)\in M,J_s\sub J;w=z_J}
f_{J_s,J}^{z,w}\ck^z_{J_s,cu}.$$
The reverse inclusion is obvious. The proposition is proved.

\subhead 2.13\endsubhead
Let $J\sub I$ and let $w\in\cw_J$. It is known that
$\dim\un\ck^w_{J,cu}=1$ if $J=J_s$ for some $s\in\Xi$ and
$\un\ck^w_{J,cu}=0$, otherwise. It follows that 2.10,2.12 can be
rewritten as follows.

$$\un\ck^w_J=\un\ck^w_{J,cu}\op
\sum_{K\prec J}\sum_{y\in{}^K_*W;y_J=w}f_{K,J}^{y,w}\un\ck^y_{K,cu}.$$
Taking direct sum over all $w\in\cw_J$, we obtain
$$\un\ck^*_J=\un\ck_J(J)\op\sum_{K\prec J}\un\ck_J(K)\tag a$$
where for any $K\sub J$ we set
$$\un\ck_J(K)=\sum_{y\in\cw_K,w\in\cw_J;y_J=w}f_{K,J}^{y,w}
\un\ck^y_{K,cu},$$
that is,
$$\un\ck_J(K)=\sum_{y\in\cw_K;y_J\in\cw_J}f_{K,J}^{y,y_J}
\un\ck^y_{K,cu}.\tag b$$

\subhead 2.14\endsubhead
Let $J\sub I$. Using 1.12(b), 2.13(a), we have
$$\align&\un\tck_J=\sum_{J'\prec J}f_{J',J}\un\ck^*_{J'}=
\sum_{J'\prec J}\sum_{K\sub J'}f_{J',J}\un\ck_{J'}(K)\\&
=\sum_{J'\prec J}\sum_{K\sub J'}\sum_{y\in\cw_K;y_{J'}\in\cw_{J'}}
f_{J',J}f_{K,J'}^{y,y_{J'}}(\un\ck^y_{K,cu}).\endalign$$
Hence using 2.13(b) we have
$$\un\tck_J\cap\un\ck^*_J=\sum_{J'\prec J}\sum_{K\sub J'}
\sum_{y\in\cw_K;y_{J'}\in\cw_{J'},y_J\in\cw_J}f_{J',J}^{y_{J'},y_J}
f_{K,J'}^{y,y_{J'}}\un\ck^y_{K,cu}.$$
Using now 1.10(a) we obtain
$$\un\tck_J\un\cap\un\ck^*_J=
\sum_{J'\prec J}\sum_{K\sub J'}\sum_{y\in\cw_K,y_J\in\cw_J}
f_{K,J}^{y,y_J}(\un\ck^y_{K,cu})
=\sum_{J'\prec J}\sum_{K\sub J'}\un\ck_J(K)$$
hence
$$\un\tck_J\cap\un\ck^*_J=\sum_{K\prec J}\un\ck_J(K).\tag a$$

\proclaim{Theorem 2.15}Let $J\sub I$. Define a linear map
$$\op_{y\in\cw_J}\un\ck^y_{J,cu}@>>>\un\ck_J/\un\tck_J\tag a$$
as the inclusion of the left hand side into $\un\ck_J$ followed by the
quotient map $\un\ck_{J_t}@>>>\un\ck_{J_t}/\un\tck_{J_t}$. This linear
map is an isomorphism.
\endproclaim
Comparing 2.14(a),2.13(a) we deduce
$$\un\ck^*_J/(\un\tck_J\cap\un\ck^*_J)=\un\ck_J(J).\tag b$$
The obvious (injective) linear map 
$$\un\ck_J^*/(\un\tck_J\cap\un\ck^*_J)@>>>\un\ck_J/\un\tck_J$$
is also surjective. (It is enough to show that
$\un\ck_J\sub\un\ck_J^*+\un\tck_J$ or that
$\un\ck_J^\spa\sub\un\tck_J$. This follows from 1.12(a).) Hence (b) can
be rewritten as $\un\ck_J/\un\tck_J=\un\ck_J(J),$ that is
$$\sum_{y\in\cw_J}\ck^y_{J,cu}=\un\ck_J/\un\tck_J.$$
The theorem follows.

\head 3. Convolution\endhead
\subhead 3.1\endsubhead
We return to the setup in 1.1. We assume that $\tG=D=G$
so that $\t=1$. Let $J\sub I$. Let 
$$\align&Z'=\{(P,P',P'',gU_P,g'U_{P'});
(P,P',P'')\in\cp_J\T\cp_J\T\cp_J,\\&
gU_P\in G/U_P,g'U_{P'}\in G/U_{P'},gPg\i=P',g'P'g'{}\i=P''\}.\endalign$$
We define
$$Z_J\T Z_J@<b_1<<Z'@>b_2>>Z_J$$
by
$$b_1(P,P',P'',gU_P,g'U_{P'})=((P,P',gU_P),(P',P'',g'U_{P'}),$$
$$b_2(P,P',P'',gU_P,g'U_{P'})=(P,P'',g'gU_P).$$
Following \cite{L05,32.5}, for $C,C'$ in $\cd_m^{UCS}(Z_J)$ we define
$C*C'=b_{2!}b_1^*(C\bxt C')$. This is again in $\cd_m^{UCS}(Z_J)$.
We say that $C*C'$ is the convolution of $C,C'$. Now convolution
induces an $\ca$-bilinear pairing $\ck_J\T\ck_J@>>>\ck_J$ and a
$\QQ(v)$-bilinear pairing $\un\ck_J\T\un\ck_J@>>>\un\ck_J$ which makes
$\un\ck_J$ into an (associative) $\QQ(v)$-algebra (without $1$ in
general).

For $J'\sub J$ we have
$$(f_{J',J}\un\ck_{J'})*\un\ck_{J}\sub f_{J',J}\un\ck_{J'},
\un\ck_J*(f_{J',J}\un\ck_{J'})\sub f_{J',J}\un\ck_{J'}.\tag a$$
(See \cite{L05, 36.11}.)

\subhead 3.2\endsubhead
As in 2.2 we fix $V\in\fC_n$ and we assume that $\tG=D=G=Sp(V)$. Let
$F:G@>>>G$ be as in 2.2. Let $t\in\Xi,t=n-(k^2+k),k\in\NN$ and let
$J_t\sub I$ be as in 2.4.

For $x,x'$ in $\un\ck_{J_t}$ we define
$$x\dot*x'=\fra{1}{2^k(v^2-1)^t\PP(v^2)}x*x'$$
$$x\ddot*x'=\fra{1}{(v^2-1)^t\PP(v^2)}x*x'$$
where 
$$\PP(u)=u^{(k^2+k)^2-k(2k-1)-(k-1)(2k-3)-\do-1}
(u+1)^{2k}(u^2+1)^{2k-1}\do(u^{2k}+1)\in\NN[u].$$

For $f:(Z_{J_t})^F@>>>\QQ$, $f':(Z_{J_t})^F@>>>\QQ$ we define
$$f*f':(Z_{J_t})^F@>>>\QQ, f\dot*f':(Z_{J_t})^F@>>>\QQ$$
by
$$\align&(f*f')(V_*,V'_*,gU(V_*))\\&=
\sum_{(\tV_*,g'U(V_*),g''U(V'_*))}f(V_*,\tV_*,g'U(V_*))
f'(\tV_*,V'_*,g''U(\tV_*))\tag a\endalign$$
where the sum is over all $\tV_*\in\ce_t^F,g'\in G^F,g''\in G^F$ such
that $g'(V_*)=\tV_*$, $g''(\tV_*)=V'_*$, $g''g'\in U(V'_*)g=gU(V_*)$
and
$$\align&(f\dot*f')(V_*,V'_*,gU(V_*))=
\fra{D_k}{\sha(Sp(V_t^\pe/V_t)^F)(q-1)^t}(f*f')(V_*,V'_*,gU(V_*))\\&=
\fra{1}{2^k(q-1)^t\PP(q)}(f*f')(V_*,V'_*,gU(V_*))\tag b\endalign$$
where $D_k$ is the dimension of the unipotent cuspidal representation
of $Sp((V_t^\pe/V_t)^F)$.
The second equality in (b) follows from \cite{L77, (8.11.1)}.

For any $w\in\cw_{J_t}$ we shall regard
$Y^{w,y}_{t,t'}:(Z_{J_t}^w)^F@>>>\QQ$
and ${}'Y^{w,y}_{t,t'}:(Z_{J_t}^w)^F@>>>\QQ$
(with $(y,t')\in M_t$) as functions $(Z_{J_t})^F@>>>\QQ$
equal to $0$ outside $(Z_{J_t}^w)^F$.
From \cite{L03, 27.9} it is known that
$$Y^{w,1}_{t,0}\dot*Y^{w',1}_{t,0}=Y^{ww',1}_{t,0}$$
if $w,w'\in\cw_{J_t},l'(ww')=l'(w)+l'(w')$,
$$(Y^{s_i,1}_{t,0}-Y^{1,1}_{t,0})\dot*(Y^{s_i,1}_{t,0}-qY^{1,1}_{t,0})=0
$$
if $i=1,2,\do,t-1$,
$$(Y^{s'_t,1}_{t,0}-Y^{1,1}_{t,0})\dot*
(Y^{s'_t,1}_{t,0}-q^{2k+1}Y^{1,1}_{t,0})=0.$$
($s_i,s'_t$ as in 2.5.) Here $l'$ is the length function of $\cw_{J_t}$.

From 3.1(a) it follows that

(c) $\un\tck_J$ is a two-sided ideal of the algebra $\un\ck_J$.
Similarly $\tck_J$ is a two-sided ideal of the algebra $\ck_J$.

\subhead 3.3\endsubhead
Let $w_1\in\cw_{J_t},w_2\in\cw_{J_t}$,
$(y_1,t'_1)\in M_t,(y_2,t'_2)\in M_t$ be such that $t'_1+t'_2>0$.
From 2.8(a) we see that at least one of
$$gr(\cl^{w_1,y_1}_{t,t'_1}),gr(\cl^{w_2,y_2}_{t,t'_2})$$
is in $\tck_{J_t}$. Since $\tck_{J_t}$ is a two-sided ideal of
$\ck_{J_t}$ (see 3.2(c)) it follows that
$$gr(\cl^{w_1,y_1}_{t,t'_1})*gr(\cl^{w_2,y_2}_{t,t'_2})$$
is an $\ca$-linear combination of elements $gr(\cl^{w,y}_{t,t'})$ with
$w\in\cw_{J_t}$, $(y,t')\in M^!_t-\{1,0\}$ and of objects in
$UCS^w_{J_t}$ with $w\in{}^{J_t}_\spa W$. From this we can deduce that

(a) ${}'Y^{w_1,y_1}_{t,t'_1}*{}'Y^{w_2,y_2}_{t,t'_2}$ is a $\QQ$-linear
combination of elements ${}'Y^{w,y}_{t,t'}$ with $w\in\cw_{J_t}$,
$(y,t')\in M^!_t-\{1,0\}$ and of characteristic functions of objects in
$UCS^w_{J_t}$ with $w\in{}^{J_t}_\spa W$.

\subhead 3.4\endsubhead
Let $w_1,w_2$ be elements on $\cw_{J_t}$ such that
$$Y^{w_1,1}_{t,0}\dot *Y^{w_2,1}_{t,0}
=a(q)Y^{z_1,1}_{t,0}+b(q)Y^{z_2,1}_{t,0}$$
for some $z_1,z_2$ in $\cw_{J_t}$; here $a(q),b(q)$ are polynomials
in $q$.
Using 2.7(c) we see that
$$\align&(2^k{}'Y^{w_1,1}_{t,0}
+\sum_{(y,t')\in M_t^!-\{1,0\}}\ti\mu_{y,t'}(q){}'Y^{w_1,y}_{t,t'})\dot*
(2^k{}'Y^{w_2,1}_{t,0}
+\sum_{(y,t')\in M_t^!-\{1,0\}}\ti\mu_{y,t'}(q){}'Y^{w_2,y}_{t,t'})\\&
=a(q)(2^k{}'Y^{z_1,1}_{t,0}
+\sum_{(y,t')\in M_t^!-\{1,0\}}\ti\mu_{y,t'}(q){}'Y^{z_1,y}_{t,t'})
+\\& b(q)(2^k{}'Y^{z_2,1}_{t,0}
+\sum_{(y,t')\in M_t^!-\{1,0\}}\ti\mu_{y,t'}(q){}'Y^{z_2,y}_{t,t'}).
\endalign$$

From this we deduce, using 3.3(a), that
$$4^k{}'Y^{w_1,1}_{t,0}\dot*{}'Y^{w_2,1}_{t,0}-
(2^ka(q){}'Y^{z_1,1}_{t,0}+2^kb(q){}'Y^{z_2,1}_{t,0})$$
is a $\QQ$-linear combination of elements ${}'Y^{w,y}_{t,t'}$ with
$w\in\cw_{J_t}$, $(y,t')\in M^!_t-\{1,0\}$ and of characteristic
functions of objects in $UCS^w_{J_t}$ with $w\in{}^{J_t}_\spa W$.
Moreover the coefficients in this linear combinations are values at
$v^2=q$ of elements in $\QQ(v^2)$.

Since this holds with $q$ replaced by $q^s$, $s=1,2,\do$, we deduce that
(with the notation $A^w_t=gr(\cl^{w,1}_{t,0})$ for $w\in\cw_{J_t}$),
$$4^kA^{w_1}_t\dot*A^{w_2}_t-(2^ka(u^2)A^{z_1}_t+2^kb(v^2)A^{z_2}_t)
\tag a$$
is a $\QQ(v^2)$-linear combination of elements $gr(\cl^{w,y}_{t,t'})$
with $w\in\cw_{J_t}$, $(y,t')\in M^!_t-\{1,0\}$ and of objects in
$UCS^w_{J_t}$ with $w\in{}^{J_t}_\spa W$; hence (a) is in
$\un\tck_{J_t}$.

We see that in the algebra $\un\ck(J_t)/\un\tck(J_t)$ (with product induced
by $\dot*$) we have
$$A^{w_1}_t\dot*A^{w_2}_t=2^{-k}(a(u^2)A^{z_1}_t+b(v^2)A^{z_2}_t).$$
Hence in the algebra $\un\ck(J_t)/\un\tck(J_t)$ (with product induced
by $\ddot*$) we have
$$A^{w_1}_t\ddot*A^{w_2}_t=a(u^2)A^{z_1}_t+b(v^2)A^{z_2}_t.$$

\proclaim{Theorem 3.5}In the algebra $\un\ck(J_t)/\un\tck(J_t)$ (with
product induced by $\ddot*$) we have
$$A^w_t\ddot*A^{w'}_t=A^{ww'}_t$$
if $w,w'\in\cw_{J_t},l'(ww')=l'(w)+l'(w')$,
$$(A^{s_i}_t-A^1_t)\ddot*(A^{s_i}_t-v^2A^1_t)=0$$
if $i=1,2,\do,t-1$,
$$(A^{s'_t}_t-A^1_t)\ddot*(A^{s'_t}_t-v^{2(2k+1)}A^1_t)=0.$$
($s_i,s'_t$ as in 2.5.) Hence the $\QQ(v)$-algebra
$\ch_t:=\un\ck(J_t)/\un\tck(J_t)$ (with product induced by $\ddot*$) is
an Iwahori-Hecke algebra over $\QQ(v)$ of type $B_t$ with parameters
$v^2,\do,v^2,v^{2(2k+1)}$ and with unit element $A^1_t$.
\endproclaim
This follows from the formulas in 3.2, 3.4.

\head 4. Complements\endhead
\subhead 4.1\endsubhead
In this section we assume that $\tG=D=G$.
We now define a variant of the varieties $Z_J$ which is more suitable
for the study of unipotent (parabolic) character sheaves.

For $J\sub I$, we set
$$\bZ_J=\{(P,P',gR_P);(P,P')\in\cp_J\T\cp_{\e(J)},
gR_P\in G/R_P,gPg\i=P'\}.$$
(Here $R_P$ is as in 0.3.) We have an obvious map $\s_J:Z_J@>>>\bZ_J$
(a torus bundle). Let $UCS(\bZ_J)$ be the set of all simple perverse
$C$ on $\bZ_J$ (up to isomorphism) such that
$\s_J^\star(C)\in UCS_J$. 
Let $\ck(\bZ_J)$ the free $\ca$-module with basis $UCS(\bZ_J)$.
From the definitions we see that
$C\m\s_J^\star(C)$ defines a bijection $UCS(\bZ_J)@>\si>>UCS_J$
and an isomorphism $\ck(\bZ_J)@>\si>>\ck_J$.
For $J'\sub J\sub I$ we define a linear map
$\bar f_{J',J}:\ck(\bZ_{J'})@>>>\ck(\bZ_J)$
by repeating the definition of $f_{J',J}$ in 1.18 with $U_P,U_Q$
replaced by $R_P,R_Q$.
Under the isomorphisms $\ck(\bZ_{J'})@>\si>>\ck_{J'}$,
$\ck(\bZ_J)@>\si>>\ck_J$, $\bar f_{J',J}$ correspons to a
nonzero multiple of $f_{J',J}$.
For $J\sub I$ we set
$\tck_{\bZ_J}=\sum_{J'\prec J}\bar f_{J',J}(\ck_{\bZ{J'}})$,
This corresponds to $\tck_J$ under the isomorphism
$\ck(\bZ_J)@>\si>>\ck_J$.

Let $J\sub I$. We define a $\QQ(v)$-bilinear pairing $C,C'\m C*C'$,
$\un\ck_{\bZ_J}\T\un\ck_{\bZ_J}@>>>\un\ck_{\bZ_J}$ 
by repeating the definition of $*$ in 3.1 with $U_P,U_{P'}$
replaced by $R_P,R_{P'}$.

\subhead 4.2\endsubhead
As in 2.2 we fix $V\in\fC_n$ and we assume that $G=Sp(V)$. Let
$F:G@>>>G$ be as in 2.2. Let $t\in\Xi,t=n-(k^2+k),k\in\NN$ and let
$J_t\sub I$ be as in 2.4.
For $w\in\cw_{J_t}$ let $\bA^w\in\un\ck(\bZ_J)$ be the element
corresponding to $A^w\in\un\ck_J$ (see 3.4). Then
$\{\bA^w;w\in\cw_{J_t}\}$ is a basis of
$\un\ck(\bZ_J)/\un\tck(\bZ_J)$.

For $x,x'$ in $\un\ck_{\bZ_J}$ we define
$$x\ddot*x'=\fra{1}{\PP(v^2)}x*x'$$
(with $\PP$ as in 3.2; note that the factor $(v^2-1)^t$ in
the definition of $x\ddot* x'$ in 3.2 is now eliminated).
This defines an algebra structure on
$\un\ck(\bZ_J)/\un\tck(\bZ_J)$.
In this algebra we have
$$\bA^w_t\ddot*\bA^{w'}_t=\bA^{ww'}_t$$
if $w,w'\in\cw_{J_t},l'(ww')=l'(w)+l'(w')$,
$$(\bA^{s_i}_t-\bA^1_t)\ddot*(\bA^{s_i}_t-v^2\bA^1_t)=0$$
if $i=1,2,\do,t-1$,
$$(\bA^{s'_t}_t-\bA^1_t)\ddot*(\bA^{s'_t}_t-v^{2(2k+1)}\bA^1_t)=0.$$
($s_i,s'_t$ as in 2.5.) This gives another realization of the
algebra $\ch_t$ in 3.5.

\subhead 4.3\endsubhead
We return to the setup in 4.1. Let $J\sub I$.
The following result extending 2.13(a) will be proved elsewhere.
$$\un\ck^*_J=\op_{K\sub J}\sum_{y\in\cw_K;y_J\in\cw_J}f_{K,J}^{y,y_J}
\un\ck^y_{K,cu}.\tag a$$
A special case of (a) (with $J=I$) is
$$\un\ck(G)=\op_{K\sub I}\sum_{y\in\cw_K}f_{K,I}^{y,1}\un\ck^y_{K,cu}.
\tag b$$
(Note that $\cw_I=W$.) In fact (a) can be deduced from (b) and its
variant where $G$ is replaced by $D$ as in 1.1.

From (a) and 1.12 one can deduce as in 2.14, 2.15 that
$$\op_{y\in\cw_J}\un\ck^y_{J,cu}@>\si>>\un\ck_J/\un\tck_J.\tag c$$
One can then deduce a  description of the algebra $\ch_q$ (see 0.1)
in terms of perverse sheaves analogous to that in 3.5.

\enddocument